\documentclass[a4paper,10pt]{article}

\usepackage{latexsym}
\usepackage{amssymb}
\usepackage{theorem}
\usepackage{amsmath}
\usepackage{amscd}
\usepackage{color}
\usepackage{graphicx}

\newtheorem{theorem}{Theorem}[section]
\newtheorem{corollary}[theorem]{Corollary}

\newtheorem{example}[theorem]{Example}
\newtheorem{proposition}[theorem]{Proposition}
\newtheorem{remark}[theorem]{Remark}
\newtheorem{definition}[theorem]{Definition}

\newcommand{\demo}{\par\noindent{\it Proof. \/}\ }
\newcommand{\enD}{\hfill $\Box$\vspace{3truemm} \par}
\newcommand{\R}{\mathbb{R}}

\newcommand{\sign}{\operatorname{sign}}
\newcommand{\bX}{\mbox{\boldmath $X$}}
\newcommand{\ba}{\mbox{\boldmath $a$}}
\newcommand{\bx}{\mbox{\boldmath $x$}}
\newcommand{\by}{\mbox{\boldmath $y$}}
\newcommand{\bz}{\mbox{\boldmath $z$}}

\newcommand{\be}{\mbox{\boldmath $e$}}
\newcommand{\bv}{\mbox{\boldmath $v$}}
\newcommand{\bw}{\mbox{\boldmath $w$}}
\newcommand{\bm}{\mbox{\boldmath $m$}}
\newcommand{\bn}{\mbox{\boldmath $n$}}

\newcommand{\bV}{\mbox{\boldmath $V$}}

\newcommand{\bU}{\mbox{\boldmath $U$}}

\begin{document}

\title{Behavior of Gauss curvatures and mean curvatures of Lightcone framed surfaces in the Lorentz-Minkowski $3$-space}

\author{Chang Xu and Liang Chen$^*$ \\
{\scriptsize{\it School of Mathematics and Statistics, Northeast Normal University, Changchun 130024, P.R.CHINA}}
}
\date{}

\maketitle

%%%%%
\begin{abstract}

In this paper, we embark on an exploration of the differential geometric properties of lightcone framed surfaces in the Lorentz-Minkowski $3$-space.	
Generally speaking, a mixed type surface is defined as a connected regular surface that features non-empty sets of spacelike and timelike points. In contrast, a lightcone framed surface is a mixed type surface that exhibits singular points, at least locally.	
To facilitate our study of the differential geometric properties of the lightcone framed surface, we introduce a valuable tool known as the modified frame along the lightcone framed surface. As a result, we manage to reveal the behavior of the Gaussian curvature and the mean curvature of the lightcone framed surface, not only at lightlike points but also at singular points.

\end{abstract}
%%%%%
{\bf Keywords:} mixed type, lightcone framed surface, curvature, lightlike point, singular point.\\

\par\noindent
{\scriptsize 2020 \textit{Mathematics Subject classification:}
	Primary 53A35; 53B30; 58K05}

\renewcommand{\thefootnote}{\fnsymbol{footnote}}
\footnote[0]{{\scriptsize $*$ Corresponding author: chenL234@nenu.edu.cn.
		
		\mbox{}\quad This work is Partially supported by  National Nature Science Foundation of China (Grant No. 12271086).}}

%%%%%%%%%%%%%%%%%%%%%%%%%%%%%%%%%%
\section{Introduction}

For a general surface in Lorentz-Minkowski $3$-space, if a point of the surface moves between spacelike and timelike regions, then there is a lightlike point. We call this kind of surface a mixed type surface that is a connected regular surface in a Lorentzian-Minkowski $3$-space with non-empty spacelike and timelike point sets \cite{Honda1, Honda2}.
The induced metric of a mixed type surface is a signature-changing metric, and their lightlike points can be regarded as singular points of this metric. So lightlike points play key roles in studying a mixed type surface. Over the past two decades, numerous studies on signature-changing metrics from different viewpoints have been conducted by several authors \cite{FRUYY,Honda,Honda3,IT,Klya,RT,Ste,Tari1,UY}. It is worth mentioning that Honda, Saji and Teramoto
employed an ingenious approach way, namely the theory of wave fronts, to investigate the properties of lightlike points by using several invariants defined by them \cite{Honda2}. Moreover, Honda, Koiso, Kokubu, Umehara and Yamada also studied the boundedness of mean curvature of mixed type surfaces and show a notable result: no constant mean curvature (CMC) mixed-type surfaces exist except for those with zero mean curvature (ZMC) \cite{Honda1}.  

On the other hand, singular points naturally appear in the mixed type surface. In this paper, we investigate the differential geometry of singular mixed type surface. In particular, we examine the properties of Gauss curvature and mean curvature at both lightlike points and singular points are considered. To this end, we first build upon the seminal work by Takahashi and his collaborators on framed surfaces in Euclidean $3$-space \cite{Fukunaga1} and the lightcone framed curve in Lorentz-Minkowski $3$-space \cite{Liang Chen1} to define the singular mixed type surface, so-called lightcone framed surface, which has both lightlike points and singular points. For details on lightcone framed surfaces in Lorentz-Minkowski $3$-space see \cite{M.Li1}. In \cite{Yao1}, Yao and Pei studied a particular kind of lightcone framed surfaces and called them generalized null framed surfaces which only consist of lightlike points and spacelike points. In this paper, however, we investigate the lightcone framed surfaces which are made up of spacelike points, timelike points and lightlike points.
Moreover, in order to investigate the differential geometric properties of the lightcone framed surface at singular points or lightlike points, we will construct a new frame, namely, the modified frame. For details, see $\S$\ref{sec3.1}. Using this modified frame, we give the definitions of Gaussian curvature and mean curvature and so on. Especially, the definitions of these curvatures are well defined even at singular points or lightlike points.

The organization of this paper is as follows.
In Section $2$, we prepare some basic notions of Lorentz-Minkowski $3$-space (cf. \cite{Liang Chen1,T.Liu1}).
Moreover, we review differential geometric properties of mixed type surfaces in Lorentz-Minkowski $3$-space developed by A. Honda, K. Saji, K. Teramoto et al. \cite{Honda1, Honda2}. Furthermore, according to our terminology, we introduce some fundamental notations on lightcone framed surfaces in Lorentz-Minkowski $3$-space developed by M. Li, D. Pei and M. Takahashi \cite{M.Li1}.
In Section $3$, we consider the local differential geometric properties of lightcone framed surfaces by using an useful tool so called modified frame. The Gaussian curvature and the mean curvature of the lightcone framed surface can be well defined at not only lightlike points but also singular points bu using this frame. 
Finally, in Section $4$, we study the properties of bounded Gaussian curvature and bounded mean curvature at singular points or lightlike points. In the final section, we give an example of the lightcone framed surface to support our results.

We shall assume throughout the whole paper that all maps and manifolds are differentiable of class $C^\infty$ and arguments of functions of surfaces are all $(u,v)$ without unless stated.

%%%%%%%%%%%%%%%%%%%%%%%%%%%%%%%%%%
\section{Preliminaries}

%%%%%%%%%%%%%%%%%%%%%%%%%%%%%%%%%%
\subsection{Lorentz-Minkowski $3$-space}

Let $\R^3=\{ (x_1,x_2,x_3)|x_i\in \R, i=1,2,3 \}$ be a $3$-dimensional vector space. 
For any vectors $\bx=(x_1,x_2,x_3),\by=(y_1,y_2,y_3)$ in $\R^3$, the {\it pseudo inner product} of $\bx$ and $\by$ is defined to be 
$\langle \bx,\by \rangle=-x_1y_1+x_2y_2+x_3y_3$. 
We call $(\R^3,\langle , \rangle)$ the {\it Lorentz-Minkowski 3-space} and write $\R_1^3$ instead of $(\R^3,\langle ,\rangle)$.

We say that non-zero vector $\bx$ in $\R_1^3$ is {\it spacelike, lightlike} or {\it timelike} if $\langle \bx,\bx \rangle>0$, $ \langle \bx,\bx \rangle=0$ or $\langle \bx,\bx \rangle<0$ repectively. We also define $\sign(\bx)=1$, $0$ or $-1$ if $\bx$ is spacelike, lightlike or timelike respectively.
The {\it norm} of the vector $\bx$ is defined by $ \| \bx \| =\sqrt{| \langle \bx,\bx\rangle|}$.

For any $\bx=(x_1,x_2,x_3)$, $\by=(y_1,y_2,y_3)\in \R_1^3$, we define a vector $\bx\wedge\by$ by 
\begin{eqnarray*}
\bx\wedge\by=\begin{vmatrix}-\be_1
&\be_2  &\be_3 \\x_1
&x_2  &x_3 \\y_1
&y_2  &y_3
\end{vmatrix},
\end{eqnarray*}
where $\{\be_1, \be_2, \be_3\}$ is the canonical basis of $\R_1^3$. 

For any $\bz \in \R_1^3$, we can easily check that 
\begin{eqnarray*}
\langle \bz,\bx\wedge\by \rangle=\det(\bz,\bx,\by),
\end{eqnarray*}
so that $\bx \wedge\by$ is pseudo-orthogonal to both $\bx$ and $\by$, and
\begin{eqnarray*}
&&(\bx\wedge\by)\wedge\bz=-\langle \bx,\bz \rangle\by+\langle \by,\bz \rangle\bx,\\
&&\bx\wedge(\by\wedge\bz)=\langle \bx,\by \rangle\bz-\langle \bx,\bz \rangle\by.
\end{eqnarray*}
Morever for any $\ba\in\R_1^3$, we have
\begin{eqnarray*}
(\bx\wedge\by,\bz\wedge\ba)=\langle \bx,\ba \rangle \langle\by,\bz\rangle-\langle \bx,\bz \rangle \langle\by,\ba\rangle.
\end{eqnarray*}

If $\bx$ is a timelike vector, $\by$ is a spacelike vector and $\bx\wedge\by=\bz$, then by a straightforward calculation, we have 
\begin{eqnarray*}
\bz\wedge\bx=\by, \by\wedge\bz=-\bx.
\end{eqnarray*}
If $\bx$ is a spacelike vector, $\by$ is a timelike vector and $\bx\wedge\by=\bz$, then by a straightforward calculation, we have 
\begin{eqnarray*}
\bz\wedge\bx=-\by, \by\wedge\bz=\bx.
\end{eqnarray*}
If both $\bx$ and $\by$ are spacelike vectors and $\bx\wedge\by=\bz$, then by a straightforward calculation, we have 
\begin{eqnarray*}
 \bz\wedge\bx=-\by, \by\wedge\bz=-\bx.
\end{eqnarray*}

For a vector $\bv \in \R_1^3$ and a real number c, we define the {\it plane} with the pseudo-normal $\bv$ by 
\begin{eqnarray*}
 P(\bv,c)=\left \{ \bx\in \mathbb{R}_1^3|\left \langle  \bx,\bv\right \rangle=c \right \}.
\end{eqnarray*}
We call $P(\bv,c)$ a {\it timelike plane}, {\it lightlike plane} or {\it spacelike plane} if $\bv$ is spacelike, lightlike or timelike, respectively.

We define {\it hyperbolic $2$-space} by
\begin{eqnarray*}
H^2(-1)=\{\bx\in \R_1^3\ |\ \langle\bx,\bx\rangle =-1\},
\end{eqnarray*}
{\it de Sitter $2$-space} by
\begin{eqnarray*}
S^2_1=\{\bx\in \R_1^3\ |\ \langle\bx,\bx\rangle =1\},
\end{eqnarray*}
{\it {\rm (}open{\rm )} lightcone} at the origin by
\begin{eqnarray*}
LC^{*}=\{\bx \in \R_1^3\setminus \{0\} \ |\ \langle \bx, \bx\rangle=0\}.
\end{eqnarray*}
In this paper, we consider a double Legendrian fibration as follows \cite{CI09,S09}:
\begin{eqnarray*}
&& \Delta_4 =\{(\bv,\bw) \in LC^* \times LC^* \ | \ \langle \bv,\bw \rangle=-2 \},\\
&& \pi_{41}:\Delta_4 \to LC^*, \pi_{41}(\bv,\bw)=\bv, \ \pi_{42}:\Delta_4 \to LC^*, \pi_{42}(\bv,\bw)=\bw, \\
&& \theta_{41}=\langle d\bv, \bw \rangle|_{\Delta_4}, \theta_{42}=\langle \bv,d\bw \rangle|_{\Delta_4}.
\end{eqnarray*}

%%%%%%%%%%%%%%%%%%%%%%%%%%%%%%%%%%
\subsection{Mixed type surfaces in Lorentz-Minkowski $3$-space}
In this subsection, we will review some basic notations of mixed type surfaces in Lorentz-Minkowsi $3$-space. For details, please refer to \cite{Honda1, Honda2}. 

Let $f: U\to \Bbb R_1^3$ be an immersion, where $U\subset \Bbb R^2$ is an open subset, we call $f$ a {\it surface} immersed in $\Bbb R_1^3$. 
For a point $p\in U$, we call it a {\it spacelike} (respectively, {\it timelike}, {\it lightlike}) point, if the first fundamental form of $f$ at point $p$ is a positive definite (respectively, indefinite, degenerate) symmetric bilinear form on $T_pU$. We denote $U_+$ (respectively, $U_-$, LD) as the set of spacelike (respectively, timelike, lightlike) points. A surface $f: U\to \Bbb R_1^3$ is called a {\it mixed type surface}, if both the sets of spacelike points and timelike points are non-empty. In other words, a mixed type surface in $\Bbb R_1^3$ is a connected regular surface whose spacelike and timelike point sets are both non-empty. 

On the other hand, although mixed type surfaces are regular smooth maps, we may regard the lightlike points of a mixed type surface as singular points of the first fundamental form as follows: we define a smooth function $\lambda: U\to \Bbb R$ by $\lambda(u,v)=E(u,v)G(u,v)-F^2(u,v)$, where d$s^2=Edu^2+2Fdudv+Gdv^2$ is the first fundamental form of $f$. 
We call $\lambda$ the {\it discriminant function}. A point $p\in U$ is called a {\it lightlike point} ({\it spacelike point or timelike point}) if and only if $\lambda(p)=0\ (\lambda(p)>0$ or $\lambda(p)<0$, respectively$)$. Moreover, a lightlike point $p$ is said to be {\it non-degendrate} if d$\lambda(p)\neq 0$. Furthermore, for a non-degenerate lightlike point $p$, by the implicit function theorem, there exists a regular curve $\gamma(t) (|t|<\epsilon)$ in $U$ such that $p=\gamma(0)$ and the image of the curve $\gamma$ consists with $LD$ on a neighborhood of $p$. $\gamma(t)$ is called a {\it characteristic curve}.

Let $f: U\to \Bbb R_1^3$ be a mixed type surface. We consider a smooth vector field $\eta$ defined on a neighborhood $W$ of $p\in LD$.
It is called a {\it null vector field} if $\eta_q$ satisfies $({\rm d}s^2)_q(\eta_q, \bv)=0$ for any $\bv\in T_qU$ and $q\in LD\cap W$. 
Moreover, for the characteristic curve $\gamma(t)$, we call $\eta(t):=\eta_{\gamma(t)}\in T_{\gamma(t)}U$ a {\it null vector field along $\gamma(t)$}. Furthermore, we can give the classification of the non-degenerate lightlike points as follows. 
\begin{definition}[\cite{ Honda2}, Definition 2.2]\label{def2.1}
	Let $p\in U$ be a non-degenerate lightlike point, $\gamma(t) (|t|<\epsilon)$ a characteristic curve passing through $p=\gamma(0)$, and $\eta(t)$ a null vector field along $\gamma(t)$. If $\gamma^{\prime}(0)$ and $\eta(0)$ are linear independent (linear dependent, respectively), we call $p$ the {\it lightlike point of the first kind} ({\it lightlike point of the second kind}, respectively).
\end{definition}

For a lightlike point $p$ of the second kind, if there exists a sequence $\{p_n\}$ of the lightlike points of the first kind such that $\lim_{n\to \infty}p_n=p$, we call $p$ the {\it admissible point}. Moreover, $p$ is said to be an $L_{\infty}$-{\it point} if $p$ is not admissible. Then we can show the following 	behavior of Gaussian curvature $K$ and mean curvature $H$ of the mixed type surface $f$. (See \cite{Honda1}, Proposition 3.5 and \cite{Honda2}, Theorem B)

\begin{theorem}
	Let $f: U\to \Bbb R_1^3$ be a mixed type surface and $p\in U$ a non-degenerate lightlike point. Then we have the following assertions:\\	
	$(1)$ If the Gaussian curvature $K$ of $f$ is bounded on a neighborhood of $p$, then $p$ must be of the first kind.\\
	$(2)$ If the mean curvature $H$ of $f$ is bounded on a neighborhood of $p$, then $p$ must be an  $L_{\infty}$-point.
\end{theorem}

%%%%%%%%%%%%%%%%%%%%%%%%%%%%%%%%%%

%%%%%%%%%%%%%%%%%%%%%%%%%%%%%%%%%%
\subsection{Lightcone framed surfaces in the Lorentz-Minkowski $3$-space}

In order to investigate mixed type of surfaces with singular points in the Lorentz-Minkowski $3$-space, we introduce the basic notations on the lightcone framed surface. 

\begin{definition}[\cite{M.Li1}]
Let $(\bX,\bv,\bw):U \to \R^3_1 \times \Delta_4$ be a smooth mapping.
We say that $(\bX,\bv,\bw):U \to \R^3_1 \times \Delta_4$ is a {\it lightcone framed surface} if there exist smooth functions: $\alpha,\beta: U\to \R$ such that 
\begin{eqnarray*}
\bX_u\wedge\bX_v=\alpha\bv+\beta\bw
\end{eqnarray*}
for all $(u,v)\in U$.
We also say that $\bX:U \to \R^3_1$ is a {\it lightcone framed base surface} if there exists a smooth mapping $(\bv,\bw):U \to \Delta_4$ such that $(\bX,\bv,\bw):U \to \R^3_1 \times \Delta_4$ is a lightcone framed surface.
\end{definition}

We denote $\bm=-(1/2)\bv \wedge \bw$.
It is obvious that $\langle\bm,\bm\rangle=1$.
We call $\{\bv,\bw,\bm\}$ the {\it lightcone frame} of $\Bbb R_1^3$ along $\bX$.
By a direct calculation, we have the following differential equations:

\begin{eqnarray}\label{curvature-one}
	\begin{pmatrix}
		\bX_u\\\bX_v
	\end{pmatrix}=\begin{pmatrix}
		a_1&  b_1& c_1\\
		a_2&  b_2& c_2
	\end{pmatrix}\begin{pmatrix}
		\bv\\
		\bw\\
		\bm
	\end{pmatrix}
\end{eqnarray}

\begin{eqnarray}\label{curvature-two}
	\begin{pmatrix}
		\bv_u\\\bw_u\\\bm_u
	\end{pmatrix}=\begin{pmatrix}
		e_1&  0& 2g_1\\
		0& -e_1& 2f_1\\
		f_1& g_1& 0
	\end{pmatrix}\begin{pmatrix}
		\bv\\
		\bw\\
		\bm
	\end{pmatrix},\begin{pmatrix}
		\bv_v\\\bw_v\\\bm_v
	\end{pmatrix}=\begin{pmatrix}
		e_2&  0& 2g_2\\
		0& -e_2& 2f_2\\
		f_2& g_2& 0
	\end{pmatrix}\begin{pmatrix}
		\bv\\
		\bw\\
		\bm
	\end{pmatrix}
\end{eqnarray}

where 
\begin{eqnarray*}
	&&a_1=-\frac{1}{2}\langle \bX_u,\bw \rangle,\ b_1=-\frac{1}{2}\langle \bX_u,\bv \rangle,\ c_1=\langle \bX_u,\bm \rangle,\\
	&&a_2=-\frac{1}{2}\langle \bX_v,\bw \rangle,\ b_2=-\frac{1}{2}\langle \bX_v,\bv \rangle,\ c_2=\langle \bX_v,\bm \rangle,\\
	&&e_1=\frac{1}{2}\langle \bv,\bw_u \rangle,\ f_1=\frac{1}{2}\langle \bw_u,\bm \rangle,\ g_1=\frac{1}{2}\langle \bv_u,\bm \rangle,\\
	&&e_2=\frac{1}{2}\langle \bv,\bw_v \rangle,\ f_2=\frac{1}{2}\langle \bw_v,\bm \rangle,\ g_2=\frac{1}{2}\langle \bv_v,\bm \rangle.
\end{eqnarray*}
Note that $a_i,b_i,c_i,e_i,f_i,g_i:U \to \R,i=1,2$ are smooth functions.
We denote 
\begin{eqnarray*}
	\mathcal G=\begin{pmatrix}
		a_1&  b_1& c_1\\
		a_2&  b_2& c_2
	\end{pmatrix},
	\mathcal F_1=\begin{pmatrix}
		e_1&  0& 2g_1\\
		0& -e_1& 2f_1\\
		f_1& g_1& 0
	\end{pmatrix},
\mathcal F_2=\begin{pmatrix}
		e_2&  0& 2g_2\\
		0& -e_2& 2f_2\\
		f_2& g_2& 0
	\end{pmatrix}.
\end{eqnarray*} 
and call $(\mathcal G,\mathcal F_1,\mathcal F_2)$ {\it basic invariants} of the lightcone framed surface $(\bX,\bv,\bw)$.
Li, Pei and Takahashi also show the Uniqueness and Existence theorems for the lightcone framed surfaces by using these invariants in \cite{M.Li1}.

On the other hand, we remark that if $(\bX,\bv,\bw)$ is a lightcone framed surface, then it satisfies one of the following two conditions at least locally:
\begin{eqnarray*}
(1) & a_1=b_1=0 \ or \ a_2=b_2=0,\\
(2) & a_1=a_2=0 \ or \ b_1=b_2=0.
\end{eqnarray*}
If $\bX$ satisfies the condition $(2)$, we can show that $\bX$ contains only lightlike points, except singular points.
Yao and Pei had studied this special kind of lightcone framed surfaces and called them generalized null framed surfaces.
They also classified them into three kinds, namely null plane, nullcone or tangent developable of a null type curve, for details see \cite{Yao1}.
In this paper, however, we investigate the lightcone framed surface which consists of spacelike points, timelike points and lightlike points together.
Therefore, without loss of generality, we may assume that $a_2=b_2\equiv0$ for all $(u,v)\in\bU$ in this paper.
Moreover, by a straightforward calculation, we can show that
$$
\bv\wedge\bm=-\bv, \bw\wedge\bm=\bw.
$$
It follows that $\bX_u\wedge\bX_v=-a_1c_2\bv+b_1c_2\bw$.
Then $p=(u,v)$ is a lightlike point of the lightcone framed surface $(\bX,\bv,\bw)$ if and only if $a_1(p)=0,b_1(p)\ne 0, c_2(p)\ne0$ or $a_1(p)\ne0,b_1(p)=0, c_2(p)\ne0$.
We denote the set of lightlike points by $L(\bX)=L_1(\bX)\cup L_2(\bX)$, where
\begin{eqnarray*}
&&L_1(\bX)=\{p\in\bU \ | \ c_2(p)\ne0,a_1(p)=0,b_1(p)\ne0 \},\\
&&L_2(\bX)=\{p\in\bU \ | \ c_2(p)\ne0,a_1(p)\ne0,b_1(p)=0 \}.
\end{eqnarray*}
Moreover, $p$ is a singular point of the lightcone framed surface $(\bX,\bv,\bw)$ if and only if $a_1^2(p)+b_1^2(p)=0$ or $c_2(p)=0$.
We denote the set of singular points by $S(\bX)=S_1(\bX)\cup S_2(\bX)$, where
\begin{eqnarray*}
&&S_1(\bX)=\{ p\in\bU \ | \ a_1^2(p)+b_1^2(p)\ne0, c_2(p)=0 \},\\
	&&S_2(\bX)=\{ p\in\bU \ | \ a_1^2(p)+b_1^2(p)=0 \}.
\end{eqnarray*}
Let $p\in\bU$ be a singular point of the lightcone framed surface $(\bX,\bv,\bw)$, if $p\in S_i(\bX)$, we call $p$ {\it the $1$-st or $2$-nd singular point} if $i=1$ or $2$, respectively.
In this paper, we only investigate the lightcone framed surface with $1$-st singular points. For $2$-nd singular points, the situation is quite different, for details, please see\cite{CX}.

%%%%%%%%%%%%%
\section{Local differential geometry of Lightcone framed surfaces}

%%%%%%%%%%%%%%%%%%%%%%%%%%%%%%%%%%
\subsection{Modified frames along lightcone framed surfaces}\label{sec3.1}

In order to investigate the differential geometric properties of the lightcone framed surfaces at lightlike points or singular points, we construct a useful frame instead of lightcone frame. 
We first define two unit normal vector fields of the lightcone framed surface $(\bX,\bv,\bw)$ with $a_2=b_2\equiv0$ by
\begin{eqnarray*}
\bn^\pm=\frac{\bX_u(u,v)\wedge\bX_v(u,v)}{||\bX_u(u,v)\wedge\bX_v(u,v)||}=\frac{-a_1(u,v)\bv(u,v)+b_1(u,v)\bw(u,v)}{||4a_1(u,v)b_1(u,v)||}
\end{eqnarray*}
at any non-lightlike point and non-singular point.
We only take $"+"$ for the convenience and denote it by $\bn$.
It follows that $\{ \bX_u,\bX_v,\bn \}$ is a frame of the lightcone framed surface at any non-lightlike point and non-singular point.
Since $\bX_v(u,v)=c_2(u,v)\bm(u,v)$ for any $(u,v)\in\bU$, we have $\bX_v(u,v)=0$ at any $1$-st singular point $(u,v)\in S_1(\bX)$.
Moreover, it is obvious that
$$
\frac{\bX_u(u,v)\wedge\bm(u,v)}{||\bX_u(u,v)\wedge\bm(u,v)||}=\frac{\bX_u(u,v)\wedge\bX_v(u,v)}{||\bX_u(u,v)\wedge\bX_v(u,v)||}
$$
and $\bX_u(u,v)\wedge\bm(u,v)=-a_1(u,v)\bv(u,v)+b_1(u,v)\bw(u,v)$ at any $(u,v)\in\bU$.
We denote 
$$
\tilde{\bn}(u,v)=\bX_u(u,v)\wedge\bm(u,v)
$$
and use $\{ \bX_u,\bm,\tilde{\bn} \}$ instead of $\{ \bX_u,\bX_v,\bn \}$ at any point $(u,v)\in\bU$ and call it {\it modified frame}.

We denote the first and second fundamental invariants by $E,F,G,L,M,N$ according to the frame $\{ \bX_u,\bX_v,\bn \}$ and $\tilde{E}, \tilde{F},\tilde{G},\tilde{L},\tilde{M},\tilde{N}$ according to the frame $\{ \bX_u,\bm,\tilde{\bn} \}$.
By a straightforward calculation, we have
\begin{eqnarray}
&&\tilde{E}=\langle \bX_u,\bX_u \rangle=c^2_1-4a_1b_1,\\
&&\tilde{F}=\langle \bX_u,\bm \rangle=c_1,\\
&&\tilde{G}=\langle \bm,\bm \rangle=1,\\
&&\tilde{L}=\langle \bX_{uu},\tilde{\bn} \rangle=2a_1(b_{1u}-b_1e_1+c_1g_1)-2b_1(a_{1u}+a_1e_1+c_1f_1),\\
&&\tilde{M}=\langle \bm_u,\tilde{\bn} \rangle=2(a_1g_1-b_1f_1),\\
&&\tilde{N}=\langle \bm_v,\tilde{\bn} \rangle=2(a_1g_2-b_1f_2).
\end{eqnarray}
We define a function $\tilde{\lambda}:=\tilde{E}\tilde{G}-\tilde{F}^2=-4a_1b_1$ for each point $(u,v)\in U$ and call it {\it modified discriminant function} of the lightcone framed surface $(\bX,\bv,\bw)$.
It follows that 
\begin{eqnarray}\label{first fundamental invariants}
E=\tilde{E}, F=c_2\tilde{F}, G=c_2^2\tilde{G},
\end{eqnarray}
\begin{eqnarray}\label{second fundamental invariants}
L=\frac{1}{\sqrt{|\tilde{\lambda}|}}\tilde{L}, M=\frac{c_2}{\sqrt{|\tilde{\lambda}|}}\tilde{M}, N=\frac{c_2}{\sqrt{|\tilde{\lambda}|}}\tilde{N}.
\end{eqnarray}
We remark that a lightlike point $p=(u,v)\in L(\bX)$ is non-degenerate if $\tilde{\lambda}(p)=0$ and ${\rm d}\tilde{\lambda}(p)\ne0$ according to  the terminology from \cite{Honda2}.

Now we consider the Gaussian curvature and the Mean curvature at any non-lightlike point and non-singular point. 
For details, please refer to \cite{Kim1}.
$$
K=\frac{LN-M^2}{|EG-F^2|},\ H=\frac{LG-2MF+NE}{2|EG-F^2|}.
$$
A direct calculation show the following:
$$
K=\frac{\tilde{L}\tilde{N}-c_2\tilde{M}^2}{c_2|\tilde{\lambda}|^2},\
H=\frac{c_2\tilde{L}\tilde{G}-2c_2\tilde{M}\tilde{F}+\tilde{N}\tilde{E}}{2c_2|\tilde{\lambda}|^\frac{3}{2}}.
$$
We define two smooth functions by
\begin{eqnarray}\label{tilde Gauss and mean}
\tilde{K}=\tilde{L}\tilde{N}-c_2\tilde{M}^2,
\tilde{H}=\frac{1}{2}(c_2\tilde{L}\tilde{G}-2c_2\tilde{M}\tilde{F}+\tilde{N}\tilde{E}).
\end{eqnarray}
We call $\tilde{K}$ {\it modified Gauss curvature} and $\tilde{H}$ {\it modified mean curvature} of the lightcone framed surfaces $(\bX,\bv,\bw)$.
By a straightforward calculation, one can show the following:
\begin{eqnarray*}
\tilde{K}=&&4a_1^2[g_2(b_{1u}-b_1e_1+c_1g_1)-c_2g_1^2]+4b_1^2[f_2(a_{1u}+a_1e_1+c_1f_1)-c_2f_1^2]\\
&&-4a_1b_1[g_2(a_{1u}+a_1e_1+c_1f_1)+f_2(b_{1u}-b_1e_1+c_1g_1)-2c_2f_1g_1],\\
\tilde{H}=&&a_1[c_2(b_{1u}-b_1e_1-c_1g_1)+g_2(c_1^2-4a_1b_1)]- b_1[c_2(a_{1u}+a_1e_1-c_1f_1)+f_2(c_1^2-4a_1b_1)].
\end{eqnarray*}

\begin{remark}
According to the above notations, modified Gauss curvature $\tilde{K}$ and modified mean curvature $\tilde{H}$ are well defined even if at lightlike points and $1$-st singular points of $\bX$.
\end{remark}

%%%%%%%%%%%%%%%%%%%%%%%%%%%%%%%%%%
\subsection{Properties of lightlike points and singular points of lightcone framed surfaces}

We first show the properties of lightlike points of the lightcone framed surface $(\bX,\bv,\bw)$ with $a_2=b_2\equiv0$ by referring to \cite{Honda2}.

As $EG-F^2=-4a_1b_1c_2^2=0$ is determined by $a_1=0$ or $b_1=0$ at lightlike point of $\bX$, we use the {\it modified discriminant function} $\tilde{\lambda}=-4a_1b_1$ to describe properties of lightlike points of $\bX$.
Then a lightlike point $p\in U$ of the lightcone framed surface $(\bX,\bv,\bw)$ is {\it non-degenerate} if d$a_1(p)\ne0$ or d$b_1(p)\ne0$, {\it denegerate} if d$a_1(p)=0$ or d$b_1(p)=0$.

\begin{proposition}
Let $p=(u,v)\in U$ be a lightlike point of the lightcone framed surface $(\bX,\bv,\bw)$.
Then the null vector at $p=(u,v)$ is 
\begin{eqnarray*}
\eta=c_2\partial_u-c_1\partial_v.
\end{eqnarray*}
\end{proposition}
%%%
\demo
Let the null vector at $p=(u,v)$ be $\eta=\eta_1\partial_u+\eta_2\partial_v$.
Then 
\begin{eqnarray*}
{\rm d}\bX(\eta)=\eta_1\bX_u+\eta_2\bX_v=\eta_1(a_1\bv+b_1\bw+c_1\bm)+\eta_2c_2\bm=\eta_1a_1\bv+\eta_1b_1\bw+(\eta_1c_1+\eta_2c_2)\bm.
\end{eqnarray*}
As the vector d$\bX(\eta)$ at $p$ is lightlike, we have
\begin{eqnarray}\label{null vector of X}
(\langle {\rm d}\bX(\eta), {\rm d}\bX(\eta) \rangle)(p)=(\eta_1^2(c_1^2-4a_1b_1)+\eta_2^2c_2^2+2\eta_1\eta_2c_1c_2)(p)=0.
\end{eqnarray}
It is obvious that $\eta_1=c_2$ and $\eta_2=-c_1$ satisfy \eqref{null vector of X}.
So we get the null vector $\eta=c_2\partial_u-c_1\partial_v$ at a lightlike point $p=(u,v)$.
\enD
%%%%
According to Definition \ref{def2.1}, we can show the following proposition.
\begin{proposition}
Let $p=(u,v)\in U$ be a non-degenerate lightlike point of the lightcone framed surface $(\bX,\bv,\bw)$. 
Then $p=(u,v)$ is the first kind lightlike point if 
$(a_{1u}c_2-a_{1v}c_1)(p)\ne0$ or $(b_{1u}c_2-b_{1v}c_1)(p)\ne0$;
the second kind lightlike point if
$(a_{1u}c_2-a_{1v}c_1)(p)=0$ or $(b_{1u}c_2-b_{1v}c_1)(p)=0$.
\end{proposition}

Next we characterize $1$-st singular points of lightcone framed surfaces $(\bX,\bv,\bw)$ with $a_2=b_2\equiv0$ by using the method of frontals in \cite{Martins1}.
A map $f:U\to\R^3_1$ is a {\it frontal} if $F=(f,[\nu])$ is a isotropic map, namely, $\langle {\rm d}f(X),\nu\rangle=0$ for any $p\in U,\ X\in T_pU$ \cite{Honda3}.
The vector field $\nu$ is the normal vector of $f$.
It is obvious that the lightcone framed base surface $\bX$ is a frontal with the normal vector $\tilde{\bn}$.

As the normal vector $\tilde{\bn}$ can be well-defined even at $1$-st singular points of $\bX$, we define the signed area density function of $\bX$ by $\Lambda:={\rm det}(\bX_u,\bX_v,\tilde{\bn})=4a_1b_1c_2$ for each point $(u,v)\in U$.
Since $\Lambda=0$ is only determined by $c_2=0$ at $1$-st singular points, we denote $\lambda_1=c_2$ to be the {\it modified signed area density function} of $\bX$ with respect to $1$-st singular points.
Then a $1$-st singular point $p\in U$ of $\bX$ is called {\it non-degenerate $1$-st singular point} $($resp. {\it degenerate $1$-st singular point}$)$ if d$c_2\ne0$ $($resp. d$c_2=0$$)$. 

Let $p$ be a non-degenerate $1$-st singular point of $\bX$.
Then there exists a regular curve $\xi(t) (|t|<\epsilon)$ in $U$ such that $p=\xi(0)$ and the image of the curve $\xi$ consists with $S_1(\bX)$ on a neighborhood of $p$. 
$\xi(t)$ is called the {\it $1$-st singular curve}.
Morever, we call the image of the $1$-st singular curve $\hat{\xi}=\bX\circ\xi$ the {\it $1$-st singular locus}.

\begin{proposition}
Let $p=(u,v)\in U$ be a $1$-st singular point of the lightcone framed surface $(\bX,\bv,\bw)$.
Then the null vector field of the lightcone framed surface $(\bX,\bv,\bw)$ at $p=(u,v)$ is $\eta=\partial_v$.
\end{proposition}
\demo
Assume that $\eta=\eta_1\partial_u+\eta_2\partial_v$, then as d$\bX(\eta)=\eta_1a_1\bv+\eta_1b_1\bw+(\eta_1c_1+\eta_2c_2)\bm$ and $c_2(p)=0$, we have d$\bX(\eta)=\eta_1(a_1\bv+b_1\bw+c_1\bm)=0$.
So we can get $\eta_1=0, \eta_2=1$.
\enD

By using the criteria of singular points of frontals introduced in \cite{IS,Martins1}, we have the following assertion.
\begin{proposition}\label{pro3.5}
Assume that $p=(u,v)\in U$ is a non-denegerate $1$-st singular point of the lightcone framed surface $(\bX,\bv,\bw)$.
Then $p$ is a singular point of the first $($resp. second $)$ kind if $c_{2v}\ne0$ $($resp. $c_{2v}=0$ $)$.
\end{proposition}

Now we consider curvatures at $1$-st singular points of lightcone framed surfaces $(\bX,\bv,\bw)$.
Let $p=(u,v)\in U$ be a non-denegerate $1$-st singular point of the first kind.
We denote $\kappa_{v}, \kappa_c$, $\kappa_\Pi$ and $\kappa_t$ the {\it limiting normal curvature}, the {\it cuspidal curvature}, the {\it product curvature} and the {\it cusp-directional torsion}, respectively as follows \cite{Martins1,Martins2}:
\begin{eqnarray*}
\kappa_v=\frac{\langle \bX_{uu},\bn \rangle}{|\bX_u|^2},
\kappa_c=\frac{|\bX_u|^\frac{3}{2}\det(\bX_u,\bX_{vv},\bX_{vvv})}{|\bX_u\wedge\bX_{vv}|^\frac{5}{2}},
\kappa_\Pi=\kappa_v\kappa_c,
\end{eqnarray*}
\begin{eqnarray*}
\kappa_t=\frac{\det(\bX_u,\bX_{vv},\bX_{uvv})}{|\bX_u\wedge\bX_{vv}|^2}-
\frac{\det(\bX_u,\bX_{vv},\bX_{uu})\langle \bX_u,\bX_{vv} \rangle}{|\bX_u|^2|\bX_u\wedge\bX_{vv}|^2}.
\end{eqnarray*}
When $\tilde{E}(p)\ne0$, we define
\begin{eqnarray*}
\tilde{\kappa}_v(p)=\frac{\tilde{L}(p)}{|\tilde{E}(p)|},
\tilde{\kappa}_c(p)=\frac{2|\tilde{E}(p)|^\frac{3}{4}\tilde{N}(p)}{|c_{2v}(p)|^\frac{1}{2}},
\tilde{\kappa}_\Pi(p)=\frac{2\tilde{L}(p)\tilde{N}(p)}{|\tilde{E}(p)|^\frac{1}{4}|c_{2v}(p)|^\frac{1}{2}},
\end{eqnarray*}
\begin{eqnarray*}
\tilde{\kappa}_t(p)=\frac{|\tilde{E}(p)|(c_{2u}(p)\tilde{N}(p)+c_{2v}(p)\tilde{M}(p))-c_{2v}(p)\tilde{F}(p)\tilde{L}(p)}{c_{2v}(p)|\tilde{E}(p)|}.
\end{eqnarray*}
It follows that
$$
\kappa_v(p)=\frac{1}{|\tilde{\lambda}(p)|^\frac{1}{2}}\tilde{\kappa}_v(p),
\kappa_c(p)=\frac{1}{|\tilde{\lambda}(p)|^\frac{5}{4}}\tilde{\kappa}_c(p),
\kappa_\Pi(p)=\frac{1}{|\tilde{\lambda}(p)|^\frac{7}{4}}\tilde{\kappa}_\Pi(p),
\kappa_t(p)=\frac{1}{|\tilde{\lambda}(p)|}\tilde{\kappa}_t(p).
$$
Then we can get the following proposition:

\begin{proposition}\label{prop3.6}
Let $p$ be a non-denegerate $1$-st singular point of the first kind of the lightcone framed surface $(\bX,\bv,\bw)$ and $\tilde{E}(p)\ne0$.
Then we have the following:\\
{\rm (1)} $\tilde{K}(p)=0$ if and only if $\tilde{\kappa}_\Pi(p)=0$.\\
{\rm (2)} $\tilde{H}(p)=0$ if and only if $\tilde{\kappa}_c(p)=0$.
\end{proposition}
\demo
According to formula \eqref{tilde Gauss and mean}, we have $\tilde{K}(p)=\tilde{L}(p)\tilde{N}(p), \tilde{H}(p)=\frac{1}{2}\tilde{N}(p)\tilde{E}(p)$. Therefore, the assertions are proved.
\enD

\begin{corollary}
	With the notations as above. Under the assumption $\tilde{E}(p)\ne0$, we have $\tilde{K}(p)=0$ if $\tilde{\kappa}_v(p)=0$.
\end{corollary}

On the other hand, let $p$ be a non-denegerate $1$-st singular point of the second kind.
We denote $\mu_c$ and $\mu_\Pi$ the {\it normalized cuspidal curvature} and the  {\it normalized product curvature} \cite{Martins1} by:
$$
\mu_c=2c_2H,\
\mu_\Pi=\kappa_v\mu_c.
$$
When $\tilde{E}(p)\ne0$, we define
$$
\tilde{\mu}_c(p)=\tilde{N}(p)\tilde{E}(p),\
\tilde{\mu}_\Pi(p)=\sign(\bX_u(p))\tilde{L}(p)\tilde{N}(p).
$$
It follows that
$$
\mu_c(p)=\frac{1}{|\tilde{\lambda}(p)|^\frac{3}{2}}\tilde{\mu}_c(p),\
\mu_\Pi(p)=\frac{1}{|\tilde{\lambda}(p)|^2}\tilde{\mu}_\Pi(p).
$$
Then we have the following proposition easily:
\begin{proposition}\label{prop3.7}
Let $p$ be a non-denegerate $1$-st singular point of the second kind of the lightcone framed surface $(\bX,\bv,\bw)$ and $\tilde{E}(p)\ne0$.
Then we have the following:\\
{\rm (1)}  $\tilde{K}(p)=0$ if and only if $\tilde{\mu}_\Pi(p)=0$.\\
{\rm (2)} $\tilde{H}(p)=0$ if and only if $\tilde{\mu}_c(p)=0$.
\end{proposition}

%%%%%%%%%%%%%
\subsection{Modified principal curvatures of lightcone framed surfaces}

We consider principal curvatures of the lightcone framed surface $(\bX,\bv,\bw)$ by using Weingarten transformation.
We set the Weingarten map $W$ at any non-lightlike points and non-singular points of $\bX$ as follows:
$$
W\begin{pmatrix}
 \bX_u\\\bX_v
\end{pmatrix}=-\begin{pmatrix}
 \bn_u\\\bn_v
\end{pmatrix}.
$$
Then multiply both sides by $\begin{pmatrix}
  \bX_u&\bX_v
\end{pmatrix}$,we get
\begin{eqnarray*}
W=&&\begin{pmatrix}
  L&M \\
  M&N
\end{pmatrix}
\begin{pmatrix}
  E&F \\
  F&G
\end{pmatrix}^{-1}
=\frac{1}{\sqrt{|\tilde{\lambda}|}}\begin{pmatrix}
  \tilde{L}&c_2\tilde{M} \\
  c_2\tilde{M}&c_2\tilde{N}
\end{pmatrix}
\begin{pmatrix}
  \tilde{E}&c_2\tilde{F} \\
  c_2\tilde{F}&c_2^2\tilde{G}
\end{pmatrix}^{-1}\\
=&&\frac{1}{\tilde{\lambda}\sqrt{|\tilde{\lambda}|}}\begin{pmatrix}
  \tilde{L}\tilde{G}-\tilde{M}\tilde{F} & \frac{\tilde{M}\tilde{E}-\tilde{L}\tilde{F}}{c_2}\\
  c_2\tilde{M}\tilde{G}-\tilde{N}\tilde{F} & \frac{\tilde{N}\tilde{E}-c_2\tilde{M}\tilde{F}}{c_2}
\end{pmatrix}.
\end{eqnarray*}
We note that the eigenvalue of the matrix $W$ is the principal curvature of the surface $\bX$.
We denote $\kappa_i(i=1,2)$ as the principal curvature of the surface $\bX$ at any non-lightlike points and non-singular points.
Then by solving the following equation
\begin{eqnarray*}
|\kappa_iI-W|=\kappa_i^2-\frac{c_2\tilde{L}\tilde{G}-2c_2\tilde{M}\tilde{F}+\tilde{N}\tilde{E}}{c_2\tilde{\lambda}\sqrt{|\tilde{\lambda}|}}\kappa_i+\frac{\tilde{L}\tilde{N}-c_2\tilde{M}^2}{c_2|\tilde{\lambda}|\tilde{\lambda}^2}=0,(i=1,2),
\end{eqnarray*}
one can obtain the principal curvatures at any non-lightlike points and non-singular points as follows:
\begin{eqnarray}\label{pricurvature}
\kappa_i=\frac{\tilde{H}\pm\sqrt{\tilde{H}^2-c_2\tilde{\lambda}\tilde{K}}}{c_2\tilde{\lambda}\sqrt{|\tilde{\lambda}|}},(i=1,2).
\end{eqnarray}
We define two functions by
\begin{eqnarray}\label{modified pricurvature}
\tilde{\kappa}_i=\frac{\tilde{K}}{\tilde{H}\pm\sqrt{\tilde{H}^2-c_2\tilde{\lambda}\tilde{K}}},(i=1,2)
\end{eqnarray}
and call them {\it modified principal curvatures} of the lightcone framed surface $(\bX,\bv,\bw)$.
It follows that 
\begin{eqnarray}\label{modified pricurvature and GM}
\tilde{\kappa}_i(u,v)=\frac{\tilde{K}(u,v)}{\tilde{H}(u,v)\pm\sqrt{\tilde{H}^2(u,v)}},(i=1,2)
\end{eqnarray}
for any points $(u,v)\in L(\bX)\cup S(\bX)$ and at lease one of them is bounded.
Without loss of generality, we assume that $\tilde{\kappa}_1(u,v)=\frac{\tilde{K}(u,v)}{2\tilde{H}(u,v)}$ is bounded and $\tilde{\kappa}_2(u,v)$ is unbounded at any points $(u,v)\in L(\bX)\cup S(\bX)$.
We remark that by the above construction, the modified principle curvatures are well defined even if at the lightlike points and singular points of $\bX$.

By Proposition \ref{prop3.6} and \ref{prop3.7}, we have the following result:
\begin{proposition}\label{prop3.9}
Assume that $p=(u,v)\in U$ is a non-degenerate $1$-st singular point of the lightcone framed surface $(\bX,\bv,\bw)$ and $E(p)\neq 0$.
Then\\
$(1)$ If $p$ is the first kind and $\tilde{\kappa}_c(p)\ne0$, we have $\tilde{\kappa}_1(p)=\sign(\bX(p))\tilde{\kappa}_v(p)$ and 
$\tilde{\kappa}_2(p)$ is unbounded.\\
$(2)$ If $p$ is the second kind and $\tilde{\mu}_c(p)\ne0$, we have 
$\tilde{\kappa}_1(p)=\sign(\bX(p))\tilde{\kappa}_v(p)$ and
$\tilde{\kappa}_2(p)$ is unbounded.
\end{proposition}
\demo
Either $p$ is the first kind and $\tilde{\kappa}_c(p)\ne0$ or $p$ is the second kind and $\tilde{\mu}_c(p)\ne0$, both mean that $\tilde{N}(p)\neq 0$ under the assumption $E(p)\neq 0$.
Moreover, as $\tilde{H}(p)=\frac{1}{2}\tilde{N}(p)\tilde{E}(p)$ and $\tilde{K}(p)=\tilde{L}(p)\tilde{N}(p)$, by Equation \ref{modified pricurvature and GM},
we know $\tilde{\kappa}_i(p)=\frac{2\tilde{L}(p)\tilde{N}(p)}{\tilde{N}(p)\tilde{E}(p)\pm |\tilde{N}(p)\tilde{E}(p)|}$.

For $\tilde{\kappa}_i(p)=\frac{2\tilde{L}(p)\tilde{N}(p)}{\tilde{N}(p)\tilde{E}(p)+|\tilde{N}(p)\tilde{E}(p)|}$, 
when $\tilde{N}(p)\tilde{E}(p)>0$, $\tilde{\kappa}_i(p)=\frac{\tilde{L}(p)}{\tilde{E}(p)}$;
when $\tilde{N}(p)\tilde{E}(p)<0$, $\tilde{\kappa}_i(p)$ is unbounded.

For $\tilde{\kappa}_i(p)=\frac{2\tilde{L}(p)\tilde{N}(p)}{\tilde{N}(p)\tilde{E}(p)-|\tilde{N}(p)\tilde{E}(p)|}$, 
when $\tilde{N}(p)\tilde{E}(p)>0$, $\tilde{\kappa}_i(p)$ is unbounded;
when $\tilde{N}(p)\tilde{E}(p)<0$, $\tilde{\kappa}_i(p)=\frac{\tilde{L}(p)}{\tilde{E}(p)}$.\\
Then we get the bounded principal curvature 
$\tilde{\kappa}_1(p)=\frac{\tilde{L}(p)}{\tilde{E}(p)}=\tilde{\kappa_v}(p)$ if $\tilde{E}(p)>0$;
$\tilde{\kappa}_1(p)=\frac{\tilde{L}(p)}{\tilde{E}(p)}=-\tilde{\kappa_v}(p)$ if $\tilde{E}(p)<0$.
\enD

Next we give principal vectors of the lightcone framed surfaces $(\bX,\bv,\bw)$ as follows:

Let $\bV_i=(V_{i1},V_{i2})(i=1,2)$ be the principal vector at non-lightlike and non-singular points of $\bX$.
Then by $(II-\kappa_i I)\bV_i=\mathbf{0}, (i=1,2)$, we have
$$
\left\{\begin{matrix}
 (L-\kappa_i E)V_{11}+(M-\kappa_i F)V_{12}=0\\
 (M-\kappa_i F)V_{11}+(N-\kappa_i G)V_{12}=0
\end{matrix}\right.,(i=1,2).
$$
then the principal vector at non-lightlike and non-singular points of $\bX$ is
$$
\bV_i=(N-\kappa_i G,-M+\kappa_i F) \ or \ (M-\kappa_i F,-L+\kappa_i E),(i=1,2).
$$
We define $\tilde{\bV}_i=(\tilde{V}_{i1},\tilde{V}_{i2})(i=1,2)$ be the principal vector of $\bX$ at $(u,v)\in U$.
Then by a straightforward calculation, we have
$$
\frac{1}{|\tilde{\lambda}|^{\frac{1}{2}}}\begin{pmatrix}
  \tilde{L}-\tilde{\kappa}_i \tilde{E}&c_2(\tilde{M}-\tilde{\kappa}_i \tilde{F}) \\
  c_2(\tilde{M}-\tilde{\kappa}_i \tilde{F})&c_2(\tilde{N}-c_2\tilde{\kappa}_i \tilde{G})
\end{pmatrix}\begin{pmatrix}
 \tilde{V}_{i1}\\\tilde{V}_{i2}
\end{pmatrix}=\begin{pmatrix}
 0\\0
\end{pmatrix},(i=1,2).
$$
The relation above is equivalent to the following two equation systems:

$$
\left\{\begin{matrix}
 (\tilde{L}-\tilde{\kappa}_1 \tilde{E})\tilde{V}_{11}+c_2(\tilde{M}-\tilde{\kappa}_1 \tilde{F})\tilde{V}_{12}=0\\
 (\tilde{M}-\tilde{\kappa}_1 \tilde{F})\tilde{V}_{11}+(\tilde{N}-c_2\tilde{\kappa}_1 \tilde{G})\tilde{V}_{12}=0
\end{matrix}\right.,
$$
and
$$
\left\{\begin{matrix}
 (c_2\tilde{L}-\tilde{\kappa} \tilde{E})\tilde{V}_{21}+c_2(c_2\tilde{M}-\tilde{\kappa} \tilde{F})\tilde{V}_{22}=0\\
 (c_2\tilde{M}-\tilde{\kappa} \tilde{F})\tilde{V}_{21}+c_2(\tilde{N}-\tilde{\kappa} \tilde{G})\tilde{V}_{22}=0
\end{matrix}\right.,
$$
where $\tilde{\kappa}=c_2\tilde{\kappa}_2$ is bounded on $U$.
Then the principal vector with respect to the bounded curvature $\tilde{\kappa}_1$ of $\bX$ is 
$$
\tilde{\bV}_1=(\tilde{N}-c_2\tilde{\kappa}_1 \tilde{G},-\tilde{M}+\tilde{\kappa}_1 \tilde{F}) \ or \ 
(c_2(\tilde{M}-\tilde{\kappa}_1 \tilde{F}), -\tilde{L}+\tilde{\kappa}_1 \tilde{E}).
$$
Note that when $\tilde{N}(u,v)\tilde{E}(u,v)\ne0$, $\tilde{\bV}_1(u,v)=(\tilde{N},\frac{\tilde{L}\tilde{F}-\tilde{M}\tilde{E}}{\tilde{E}})(u,v)$ or $(0,0)$ for any $(u,v)\in S_1(\bX)$.
So we take the principal vector 
$$
\tilde{\bV}_1=(\tilde{N}-c_2\tilde{\kappa}_1 \tilde{G},-\tilde{M}+\tilde{\kappa}_1 \tilde{F}).
$$
And the principal vector with respect to the unbounded curvature $\tilde{\kappa}_2$ of $\bX$ is
$$
\tilde{\bV}_2=(c_2(\tilde{N}-\tilde{\kappa} \tilde{G}),-c_2\tilde{M}+\tilde{\kappa} \tilde{F}) \ {\rm or} \ 
(c_2(c_2\tilde{M}-\tilde{\kappa} \tilde{F}),-c_2\tilde{L}+\tilde{\kappa} \tilde{E}),
$$
where $\tilde{\kappa}=c_2\tilde{\kappa}_2$.

\begin{proposition}\label{prop3.10}
Let $p=(u,v)\in U$ be a non-degenerate $1$-st singular point of the lightcone framed surface $(\bX,\bv,\bw)$ and $\xi$ be the $1$-st singular curve passing through $p$. Suppose that $\tilde{E}(p)\ne0$, then we have the following:\\
$(1)$ If $p$ is the first kind, then $\hat{\xi}$ is a line of curvature of $\bX$ if and only if 
$$
\tilde{\kappa}_t(p)=\frac{2(\tilde{M}(p)\tilde{E}(p)-\tilde{L}(p)\tilde{F}(p))}{\tilde{E}(p)} \ or \ 
2\tilde{M}(p),
$$
$(2)$ If $p$ is the second kind, then $\hat{\xi}$ is a line of curvature of $\bX$ if and only if $\tilde{\mu}_c(p)=0$.
\end{proposition}
\demo
$(1)$ Since $p$ is the first kind, we have $c_{2v}(p)\neq 0$. When $\tilde{E}(p)>0$,
\begin{eqnarray*}
\tilde{\kappa}_t(p)=\frac{c_{2u}(p)\tilde{N}(p)\tilde{E}(p)+c_{2v}(p)(\tilde{M}(p)\tilde{E}(p)-\tilde{L}(p)\tilde{F}(p))}{c_{2v}(p)\tilde{E}(p)}
=\frac{c_{2u}(p)\tilde{V}_{11}(p)-c_{2v}(p)\tilde{V}_{12}(p)}{c_{2v}(p)}.
\end{eqnarray*}
When $\tilde{E}(p)<0$,
\begin{eqnarray*}
\tilde{\kappa}_t(p)=\frac{c_{2u}(p)\tilde{N}(p)\tilde{E}(p)+c_{2v}(p)(\tilde{M}(p)\tilde{E}(p)+\tilde{L}(p)\tilde{F}(p))}{c_{2v}(p)\tilde{E}(p)}
=\frac{c_{2u}(p)\tilde{V}_{11}(p)}{c_{2v}(p)}+\frac{c_{2v}(p)(\tilde{M}(p)\tilde{E}(p)+\tilde{L}(p)\tilde{F}(p))}{c_{2v}(p)\tilde{E}(p)}.
\end{eqnarray*}
As $\hat{\xi}$ is a line of curvature of $\bX$ if and only if $(\tilde{\bV}_1c_2)(p)=0$ if and only if $c_{2u}(p)\tilde{V}_{11}(p)=-c_{2v}(p)\tilde{V}_{12}(p)$. According to the definition of $\tilde{\kappa}_{v}$ and Proposition \ref{prop3.9}, we have $c_{2u}(p)\tilde{N}(p)=c_{2v}(p)(\tilde{M}(p)-\frac{\tilde{L}(p)}{\tilde{E}(p)}\tilde{F}(p))$. 
This means that 
$\tilde{\kappa}_t(p)=\frac{2(\tilde{M}(p)\tilde{E}(p)-\tilde{L}(p)\tilde{F}(p))}{\tilde{E}(p)}$ under the condition $\tilde{E}(p)>0$ or 
$\tilde{\kappa}_t(p)=2\tilde{M}(p)$ under the condition $\tilde{E}(p)<0$. Thus we get the result.\\
$(2)$ By Proposition \ref{pro3.5}, if $p$ is a non-degenerate $1$-st singular points of the second kind, we have $c_{2v}=0$ with $c_{2u}\ne0$.
Then $(\tilde{\bV}_1c_2)(p)=c_{2u}(p)\tilde{V}_{11}(p)=c_{2u}(p)\tilde{N}(p)=0$ if and only if $\tilde{N}(p)=0$.
So we get the assertion.
\enD

%%%%%%%%%%%%%%%%%%%%%%%%%%%%%%%%%%
\section{Behavior of Gauss curvatures and mean curvatures of lightcone framed surfaces}

In this section, we investigate the differential geometric properties of Gauss curvatures and mean curvatures of lightcone framed surfaces.

As $K=\tilde{K}/c_2|\tilde{\lambda}|^2,\ H=\tilde{H}/c_2|\tilde{\lambda}|^\frac{3}{2}$, if $p=(u,v)$ is a lightlike point of $\bX$ and $\tilde{K}(p)\ne0$ (or $\tilde{H}(p)\ne0$), then $K$ (or $H$) is unbounded at $p$.
The case for a singular point $q=(u,v)\in U$ is similar to the case for the lightlike point $p$.
Note that $K$ and $H$ cannot be both bounded on a neighborhood of a non-denegerate lightlike point $p\in U$ \cite{Honda2}. Moreover, if $H$ is bounded on a neighborhood of a non-degenerate the first kind singular point of a frontal, then $K$ is bounded $($\cite{Martins1}, Corollary $3.13)$. So that we have the following proposition.

\begin{proposition}
	Assume that $\bX: U\to \Bbb R_1^3$ is a lightcone framed surface with $a_2=b_2\equiv0$, $K$ and $H$ are the Gauss curvature and the mean curvature of $\bX$. Let $q\in S_1(\bX)$ be the $1$-st singular point of $\bX$. If $c_{2v}(q)\neq 0$ and $H$ is bounded at $q$, then  $K$ is also bounded at $q$.
\end{proposition}
%%%%%%%
\demo According to Corollary 3.13 in \cite{Martins1} and Proposition \ref{pro3.5}, we can get the proof.
\enD

The following theorem characterizing the properties of bounded mean curvatures of mixed type surfaces at lightlike points is well known. 

\begin{theorem}(cf. \cite{Honda1})
Let $f:U \to \Bbb R_1^3$ be a mixed type surface in Lorentz-Minkowski $3$-space.
Denote by $U_+$ $($resp. $U_-$$)$ the set of points where $f$ is spacelike $($resp. timelike$)$.
Suppose that $U_+,U_-$ are both non-empty and the Mean curvature $H$ on $U_+\cup U_-$ is bounded.
Then for each $p \in \overline{U}_+\cap \overline{U}_-$, there exists a sequence $\{p_n\}_{n=1,2,3,...}$ in $U_+$ $($resp. $U_-$$)$ converging to $p$ such that $\lim_{n \to \infty} H(p_n)=0$, where $\overline{U}_+,\overline{U}_-$ are the closures of $U_+,U_-$ in $U$.
\end{theorem}

We now study the properties of bounded Gauss curvatures of lightcone framed surfaces at lightlike points by using the similar arguments to those of the mean curvatures of mixed type surfaces. With the same notations as above, we have the following theorem.

\begin{theorem}\label{the Gaussian curvature of X}
	Let $(\bX,\bv,\bw):U \to \R_1^3\times\Delta_4$ be a lightcone framed surface with $a_2=b_2\equiv0$.
	Suppose that the Gaussian curvature $K$ of $\bX$ is bounded on $U_+\cup U_-$, where $U_+ (resp. \ U_-)$ is the set of spacelike $($resp. timelike$)$ points of $\bX$ and non-empty.
	Then for each $p\in L(\bX)$, there exists a sequence $\{p_n\}_{n=1,2,3,...}$ in $U_+$$($resp. $U_-$$)$ converging to $p$ such that 
	$$
	\lim_{n \to \infty} K(p_n)=0 \ or \ \lim_{n \to \infty} c_2(p_n)K(p_n)=p_0\ne0,
	$$
	where $p_0$ is a non-zero real number.
\end{theorem}
%%%%%
\demo
Let $q_1,q_2\in V\subset U$ satisfy ${\tilde{\lambda}}(q_1)>0,{\tilde{\lambda}}(q_2)<0$ and $\gamma(s),s\in[0,2\pi]$ be a smooth curve in $U$ satisfied $\gamma(0)=q_1,\gamma(2\pi)=q_2$.
Then by Fourier series expansion of $\gamma(s)$, there exists a sufficiently large positive integer $N$, such that the real analytic curve defined by $\gamma_N(s)$ satisfies
\begin{eqnarray}\label{lambda-gamma-2}
	\tilde{\lambda}(\gamma_N(0))>0,\tilde{\lambda}(\gamma_N(2\pi))<0.
\end{eqnarray}
We set 
$$
\check{\lambda}(s):=\tilde{\lambda}(\gamma_N(s)),s\in[0,2\pi]
$$
Since $\check{\lambda}(s)$ is a real analytic function defined on $[0,2\pi]$, the set of zeros of the function $\check{\lambda}(s)$ consists of a finite set of points
$$
0<s_1<s_2<\cdot\cdot\cdot<s_n<2\pi.
$$
By \eqref{lambda-gamma-2}, we can choose the number $j$ such that the sign of $\check{\lambda}(s)$ changes from positive to negative at $s=s_j$.
So there exists a positive integer $m$ such that
$$
\lim_{s \to s_j}\frac{\check{\lambda}(s)}{(s-s_j)^m}=\lambda_0\ne0,
$$
where $\lambda_0$ is a non-zero real number.
As $\tilde{K}$ can be well defined at lightlike points and $1$-st singular points of $\bX$,
we set
$$
\check{K}(s):=\tilde{K}(\gamma_N(s)),s\in[0,2\pi],
$$
then we have
$$
K(\gamma_N(s)):=\frac{\check{K}(s)}{c_2(s)|\check{\lambda}(s)|^2},s\ne s_1,s_2,\cdot\cdot\cdot,s_n
$$
where $c_2(s)\ne0$.
As $K(p)$ is bounded at each $p\in L(\bX)$, we get $\tilde{K}(s_j)=0$.
Since $\tilde{K}(s)=0$ is a real analytic function, there exists a positive integer $l$ such that 
$$
\lim_{s \to s_j} \frac{\check{K}(s)}{(s-s_j)^l}=K_0\ne0,
$$
where $K_0$ is a non-zero real number.
Then it holds that
\begin{eqnarray*}
	K(\gamma_N(s))=\frac{\check{K}(s)}{(s-s_j)^l} \frac{(s-s_j)^{2m}}{|\check{\lambda}(s)|^{2m}}\frac{(s-s_j)^{l-2m}}{c_2(s)},
\end{eqnarray*}
\begin{eqnarray*}
	\lim_{s \to s_j} c_2(s)|s-s_j|^{2m-l}|K(\gamma_N(s))|=\frac{|K_0|}{|\lambda_0|^2}\ne0.
\end{eqnarray*}
As $K(p)$ is bounded, we have $l-2m\ge0$.
If $l-2m>0$, we have 
$$
\lim_{s \to s_j} |K(\gamma_N(s))|=0;
$$
if $l-2m=0$, we have 
$$
\lim_{s \to s_j}c_2(s)|K(\gamma_N(s))|=\frac{|K_0|}{|\lambda_0|^2}\ne0.
$$
Now we set $p_n:=\gamma_N(s)$, so we obtain the result.
\enD
%%%%%
On the other hand, we also investigate the properties of bounded mean curvatures or bounded Gauss curvatures of lightcone framed surfaces at $1$-st singular points respectively. We first give the characterization of bounded mean curvatures as follows.

\begin{theorem}\label{corollary4.10} 
Let $(\bX,\bv,\bw):U \to \R_1^3\times\Delta_4$ be a lightcone framed surface with $a_2=b_2\equiv0$ and the mean curvature $H$ of $\bX$ is bounded on $U$.
Assume that $q\in S_1(\bX)$ is a $1$-st singular point of $\bX$ and $V \in U$ is a neighborhood of $q$.
Then there exists a sequence $\{q_n\}_{n=1,2,3...}$ in $V$ converging to $q$, such that
$$
\lim_{n \to \infty} H(q_n)=0 \ or \ \lim_{n \to \infty} H(q_n)=H_0\ne0.
$$
where $H_0$ is a non-zero real number.
\end{theorem}
%%%%%
\demo
As $H=\tilde{H} / c_2|\tilde{\lambda}|^{\frac{3}{2}}$, we set $\Gamma=c_2|\tilde{\lambda}|^{\frac{3}{2}}$.
It is obvious that $\Gamma(q)=0$.
So there exists a positive integer $m$ such that
$$
\lim_{n \to \infty} \frac{\Gamma(q_n)}{(q_n-q)^m}=b\ne0,
$$
where $b$ is a non-zero real number.
As $H$ is bounded on $U$ and $\Gamma(q)=0$, we have $\tilde{H}(q)=0$.
Then there exists a positive integer $l$ such that 
$$
\lim_{n \to \infty} \frac{\tilde{H}(q_n)}{(q_n-q)^l}=a\ne0,
$$
where $a$ is a non-zero real number.
Then it holds that
\begin{eqnarray*}
\lim_{n \to \infty} H(q_n)=\lim_{q_n \to q} \frac{\tilde{H}(q_n)}{(q_n-q)^l} \frac{(q_n-q)^{m}}{\Gamma(q_n)}(q_n-q)^{l-m},
\end{eqnarray*}
\begin{eqnarray*}
\lim_{n \to \infty} (q_n-q)^{m-l}H(q_n)=a/b\ne0.
\end{eqnarray*}
As $H(q)$ is bounded, we have $l-m\ge0$.
If $l-m>0$, we have 
$$
\lim_{n \to \infty} H(q_n)=0;
$$
if $l-m=0$, we have 
$$
\lim_{n \to \infty} H(q_n)=H_0\ne0.
$$
Therefore, we obtain the result.
\enD

 Moreover, we consider the properties of bounded Gauss curvatures of lightcone framed surfaces at $1$-st singular points as follows.
\begin{theorem}\label{theorem4.11}
Let $(\bX,\bv,\bw):U \to \R_1^3\times\Delta_4$ be a lightcone framed surface with $a_2=b_2\equiv0$ and the Gaussian curvature $K$ of $\bX$ is bounded on $U$.
Suppose that $q\in S_1(\bX)$ is the $1$-st singular points of $\bX$ and $V \in U$ is a neighborhood of $q$.
Then there exists a sequence $\{q_n\}_{n=1,2,3...}$ in $V$ converging to $q$, such that
$$
\lim_{n \to \infty} K(q_n)=0 \ or \ \lim_{n \to \infty} K(q_n)=K_0\ne0.
$$
where $K_0$ is a non-zero real number.
\end{theorem}
%%%%%
\demo
As $K=\tilde{K} / c_2|\tilde{\lambda}|^2$, we set $\Gamma=c_2|\tilde{\lambda}|^2$.
It is obvious that $\Gamma(q)=0$.
So there exists a positive integer $m$ such that
$$
\lim_{n \to \infty} \frac{\Gamma(q_n)}{(q_n-q)^m}=b\ne0,
$$
where $b$ is a non-zero real number.
As $K$ is bounded on $U$ and $\Gamma(q)=0$, we have $\tilde{K}(q)=0$.
Then there exists a positive integer $l$ such that 
$$
\lim_{n \to \infty} \frac{\tilde{K}(q_n)}{(q_n-q)^l}=a\ne0,
$$
where $a$ is a non-zero real number.
Then it holds that
\begin{eqnarray*}
\lim_{n \to \infty} K(q_n)=\lim_{q_n \to q} \frac{\tilde{K}(q_n)}{(q_n-q)^l} \frac{(q_n-q)^{m}}{\Gamma(q_n)}(q_n-q)^{l-m},
\end{eqnarray*}
\begin{eqnarray*}
\lim_{n \to \infty} (q_n-q)^{m-l}K(q_n)=a/b\ne0.
\end{eqnarray*}
As $K(q)$ is bounded, we have $l-m\ge0$.
If $l-m>0$, we have 
$$
\lim_{q_n \to q} K(q_n)=0;
$$
if $l-m=0$, we have 
$$
\lim_{n \to \infty} K(q_n)=K_0\ne0.
$$
So that we complete the proof.
\enD
%%%%%

\section{Example}
In this section, we give an example to illustrate our results. This example appeared in \cite{Honda2}. However, we modify it slightly according to our terminology.
\begin{example}
We define $(\bX,\bv,\bw):[-(\pi/2),\pi/2]\times[0,2\pi) \to \R^3_1 \times \Delta_4$ by
\begin{eqnarray*}
&&\bX=(\sin u, \cos u \sin v, \cos u \cos v),\\
&&\bv=(1, \sin v, \cos v),\\
&&\bw=(1,-\sin v,-\cos v).
\end{eqnarray*}
Then $(\bX,\bv,\bw):U \to \R^3_1 \times \Delta_4$ is a lightcone framed surface with lightcone frame. 
By a direct calculation, we have
\begin{eqnarray*}
\bm=(0,-\cos v,\sin v).
\end{eqnarray*}
It follows that $\{\bv,\bw,\bm\}$ is the lightcone frame of $\bX$.
The basic invariants $(\mathcal G,\mathcal F_1,\mathcal F_2)$ of $\bX$ are as follows:
\begin{eqnarray*}
\mathcal G=\begin{pmatrix}
		a_1&  b_1& c_1\\
		a_2&  b_2& c_2
	\end{pmatrix}=\begin{pmatrix}
		-\frac{1}{2}(\sin u-\cos u)&  \frac{1}{2}(\sin u+\cos u)& 0\\
		0&  0& -\cos u
	\end{pmatrix},
\end{eqnarray*}
\begin{eqnarray*}
\mathcal F_1=\begin{pmatrix}
		e_1&  0& 2g_1\\
		0& -e_1& 2f_1\\
		f_1& g_1& 0
	\end{pmatrix}=\begin{pmatrix}
		0&  0& 0\\
		0& 0& 0\\
		0& 0& 0
	\end{pmatrix},\ 
\mathcal F_2=\begin{pmatrix}
		e_2&  0& 2g_2\\
		0& -e_2& 2f_2\\
		f_2& g_2& 0
	\end{pmatrix}=\begin{pmatrix}
		0&  0& -1\\
		0& 0& 1\\
		\frac{1}{2}& -\frac{1}{2}& 0
	\end{pmatrix}.
\end{eqnarray*}
The lightlike point set of $\bX$ is
$$
L(\bX)=\{ (u,v)\in [-\frac{\pi}{2},\frac{\pi}{2}]\times[0,2\pi) \ | \ u=\pm\frac{\pi}{4} \}
$$
and the singular point set of $\bX$ is
$$
S(\bX)=\{ (u,v)\in [-\frac{\pi}{2},\frac{\pi}{2}]\times[0,2\pi) \ | \ u=\pm\frac{\pi}{2} \}.
$$
Since $a_1(u,v)=-\frac{1}{2}(\sin u-\cos u)$, $b_1(u,v)=\frac{1}{2}(\sin u+\cos u)$ and $c_2(u,v)=-\cos u$, it is clearly that a point with $u=\pm\frac{\pi}{2}$ is the $1$-st singular point of $\bX$ and there is no the $2$-nd singular points in $\bX$. See Figure 1.

Moreover, by a direct calculation, we have:
\begin{eqnarray*}
&&\tilde{\lambda}=-\cos 2u,\\
&&\tilde{\bn}=(\sin u,-\cos u\sin v,-\cos u\cos v),\\
&&\tilde{K}=-\cos u,\\
&&K=\sec^22u,\\
&&\tilde{H}=-\sin^2u\cos u,\\
&&H=\sin^2u(|\sec2u|)^{3/2}.
\end{eqnarray*}
It follows that 
$$K(\pm \frac{\pi}{2}, v)=1,\ H(\pm \frac{\pi}{2}, v)=1,\ \tilde{K}(\pm \frac{\pi}{2}, v)=0,\ \tilde{H}(\pm \frac{\pi}{2}, v)=0\ {\rm and}\ \tilde{\kappa}_1(\pm \frac{\pi}{2}, v)=\frac{1}{2}.$$
 Moreover, both of the Gauss curvature $K$ and the mean curvature $H$ are unbounded at the lightlike points of $\bX$.

As $c_{2u}(\pm\frac{\pi}{2},v)=\pm1, c_{2v}(\pm\frac{\pi}{2},v)=0$, points $(\pm\frac{\pi}{2},v)$ are all non-degenerate $1$-st singular points of the second kind of $\bX$.
By a straightforward calculation, we have
$$
\tilde{\kappa}_v(\pm\frac{\pi}{2},v)=1,\
\tilde{\mu}_c(\pm\frac{\pi}{2},v)=0,\
\tilde{\mu}_{\Pi}(\pm\frac{\pi}{2},v)=0.
$$

\begin{center}
 \centering
 \begin{minipage}[c]{0.35\textwidth}
  \centering
  \includegraphics[scale=0.5]{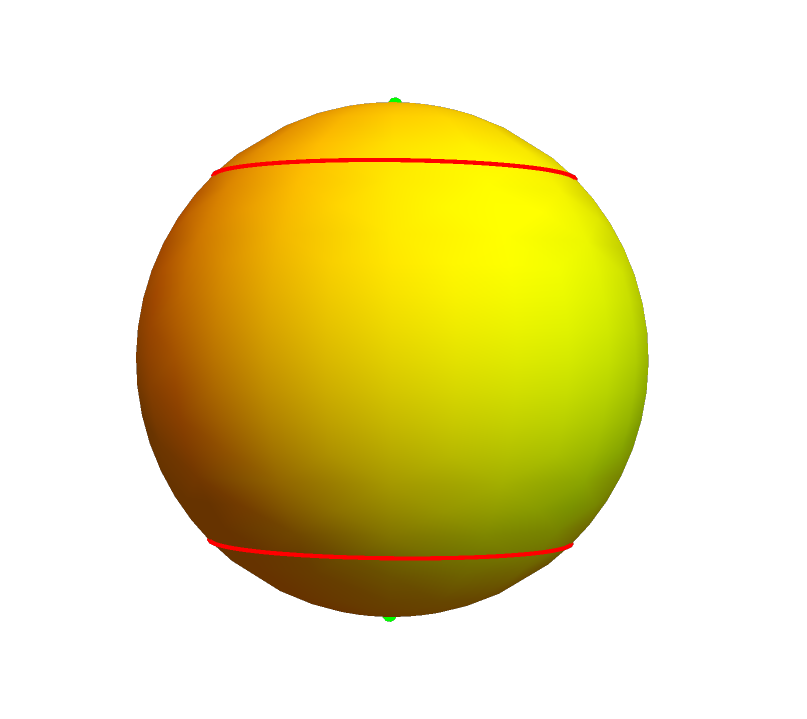} \\
 \end{minipage}
\\ Figure 1.\quad the lightcone framed surface\\ {\color{red} the lightlike locus} $($red$)$ and {\color{green} the $1$-st singular locus} $($green$)$
\end{center}

\end{example}

%%%%%%%%%%%%%%%%%%%%%%%%%%%%%%%%%

\end{document}